\documentclass[10pt]{article}%
\usepackage{amsmath}
\numberwithin{equation}{section}  
\usepackage{amsfonts}
\usepackage{amssymb}
\usepackage{amsthm}
\usepackage{graphicx}%
\providecommand{\U}[1]{\protect\rule{.1in}{.1in}}
\providecommand{\U}[1]{\protect\rule{.1in}{.1in}}
\begin{document}

\begin{center}
\thispagestyle{empty}  
\setcounter{page}{0}  

${}$
\vskip.1in
{\Large \textbf{Characteristic function and operator approach to
M-indeterminate probability densities}}

\bigskip
\bigskip

\textbf{Patrick Loughlin$^{a*}$, Leon Cohen$^{b}$}

\bigskip

\bigskip
\bigskip

$^{a}$ Departments of Electrical \& Computer Engineering, and Bioengineering, 
302 Benedum Hall, \\ University of Pittsburgh, Pittsburgh, PA 15261, USA.
Email: loughlin@pitt.edu\\
\bigskip
\bigskip
$^{b}$Department
of Physics, Hunter College of the City University of New York, 695 Park Ave.,\\
New York, NY 10065, USA. Email: leon.cohen@hunter.cuny.edu\\

\bigskip
\bigskip
\bigskip

\noindent {*Corresponding author  }

\bigskip
\bigskip
\bigskip

\noindent Declarations of interest: none

\end{center}

\newpage

\begin{center}
{\Large \textbf{Characteristic function and operator approach to
M-indeterminate probability densities}}

\bigskip

\textbf{Patrick Loughlin, Leon Cohen}

\end{center}

\bigskip

\bigskip

\noindent\textbf{{Abstract:}} Based on a quantum mechanical approach, we investigate moment- (or M-) indeterminate probability
densities by way of the characteristic function and self-adjoint operators.
The approach leads to new methods to construct classes of M-indeterminate probability
densities.

\bigskip


\section{Introduction}
In a well-known series of papers, Aharonov {\it et al.} \cite{ahar1,ahar2, ahar-book} 
developed a quantum
mechanical state function that resulted in a probability density that
depends on a parameter while the moments do not.  An analogous situation arises
in probability theory.
Specifically, while many of the common
probability densities are uniquely determined by their moments, some
are not; such densities are said to be ``moment-indeterminate,'' or
``M-indeterminate'' \cite{feller,lukacs,stoyanov2013}.  Yet, the approach of Aharonov {\it et al.} and
the way in which probability densities are obtained in quantum mechanics
is very different compared to standard probability theory.  Nevertheless, one
does not need to understand quantum theory in order to apply 
that approach to the problem of generating
M-indeterminate densities in standard probability, as we do here. In particular, we show that
using the characteristic function as defined in probability theory together
with operators that are standard in quantum mechanics allows one
to readily generate an infinite number of M-indeterminate densities.

Consider the generally complex function $g(x)$, which we shall take to be normalizable to one, 
\begin{equation}
\int_{-\infty}^{\infty}\left\vert g(x)\right\vert ^{2}dx=1 \nonumber
\end{equation}
Although not essential to our developments here, in quantum mechanical language $g(x)$ is called a state function and $x$ may be
the position of a particle.  What is important for our considerations is that 
the probability that the particle resides within some interval $x_1 < x <x_2$  is given by
\begin{equation}
\text{Pr}\{x_1 < x <x_2\}= \int_{x_1}^{x_2}\left\vert g(x)\right\vert ^{2}\,dx \nonumber
\end{equation}
Hence, the probability density is $f_X(x)=\left|g(x)\right|^2$.

To obtain the probability density of other variables,  $g(x)$ is expanded
 in a complete set of functions $u(r,x)$ \cite{morse},
\begin{equation}
g(x)=\int_{-\infty}^{\infty}G(r)u(r,x)\,dr  \label{eigen-expansion}%
\end{equation}
where
\begin{equation}
G(r)=\int_{-\infty}^{\infty}g(x)u^{\ast}(r,x)\,dx   \label{r-transform}%
\end{equation}
The functions $u(r,x)$ are  eigenfunctions that are obtained by solving the
eigenvalue problem
\begin{equation}
\mathcal{A}\,u(r,x)=r\,u(r,x)\nonumber\label{eigvalprob}%
\end{equation}
where the operator $\mathcal{A}\,\ $is self-adjoint. (In writing these
equations, we have assumed the eigenvalues are continuous.) The function $G(r)$
is the representation of $g(x)$ in the $r$-domain. Note that because $g(x)$ is
normalized, so is $G(r)$, namely
\begin{equation}
\int_{-\infty}^{\infty}\left\vert G(r)\right\vert ^{2}dr=\int_{-\infty
}^{\infty}\left\vert g(x)\right\vert ^{2}dx=1\nonumber
\end{equation}
The important point is that the probability density for $r$ is 
$f_R(r)=\left\vert G(r)\right\vert ^{2}$  \cite{bohm, merz, wilcox}.
(We also note that the transformation given by Eq. \eqref{r-transform} will always result in a
proper probability density as long as the operator is self-adjoint,
which also implies that the eigenvalues are real and the eigenfunctions are
complete and orthogonal.)

As a concrete example, consider the operator 
\begin{equation}
\mathcal{A}=\frac{1}{i}\frac{d}{dx}\nonumber
\end{equation}
which in quantum mechanics is called the momentum operator.  For this case, one has
\begin{equation}
G(r)=\frac{1}{\sqrt{2\pi }}\int_{-\infty}^\infty g(x)e^{-ixr}dx \nonumber
\end{equation}
Hence, we see that $g(x)$ and $G(r)$ are a Fourier transform pair.
Now consider the moments
\begin{equation}
\text{E\negthinspace}\left[  X^{n}\right]  \,=\,\int_{-\infty}^{\infty}%
x^{n}\left\vert g(x)\right\vert ^{2}dx\,\nonumber
\end{equation}
By virtue of the Fourier theory, these may be calculated from $G$ by way of
\begin{equation}
\text{E\negthinspace}\left[  X^{n}\right]  \,=\,\int_{-\infty}^{\infty}%
G^{\ast}(r)\left(  -\frac{1}{i}\frac{d}{dr}\right)  ^{n}G(r)\,dr\,\nonumber
\end{equation}
Similarly, one has
\begin{align}
\text{E\negthinspace}\left[  R^{n}\right] &=\int_{-\infty}^{\infty}%
r^{n}\left\vert G(r)\right\vert ^{2}dr\,\nonumber \\
&=\int_{-\infty}^{\infty
}g^\ast(x)\left(  \frac{1}{i}\frac{d}{dx}\right)  ^{n}g(x)\,dx\,\nonumber
\end{align}
We build on these relationships and their generalization for other self-adjoint operators to develop
a characteristic function approach to M-indeterminate densities.  Before that, however, we briefly discuss the function $g(x)$ considered by Aharonov and colleagues \cite{ahar1,ahar2, ahar-book}.

\subsection{The approach of Aharonov {\it et al.}}
Consider the function
\begin{equation}
g(x)=\frac{1}{\sqrt{2}}\left(  g_{1}(x)+e^{i\beta}g_{2}(x)\right)   \nonumber\label{eq1}%
\end{equation}
where $\beta$ is a real parameter and where each of $g_{1}(x)$ and $g_{2}(x)$ is
normalized to one. In addition, $g_{1}(x)$ and $g_{2}(x)$ are taken to
be finite extent and, importantly, to have disjoint support such that $g_{1}(x)g_{2}(x)=0$.
The probability density is then
\begin{equation}
f_X(x)=\left\vert g(x)\right\vert ^{2}=\frac{1}{2}\left\vert g_{1}(x)+e^{i\beta
}g_{2}(x)\right\vert ^{2}\,=\,\frac{1}{2}\left(  \left\vert g_{1}%
(x)\right\vert ^{2}+\left\vert g_{2}(x)\right\vert ^{2}\right) \nonumber
\end{equation}
and the moments are given by
\begin{equation}
\text{E\negthinspace}\left[   X^{n}\right] =\frac{1}{2}\,\int_{-\infty}^\infty x^{n}\left(  \left\vert
g_{1}(x)\right\vert ^{2}+\left\vert g_{2}(x)\right\vert ^{2}\right)\,dx \nonumber
\end{equation}
Since the probability density is independent of $\beta$, 
naturally so are the moments.

Now consider the $r$-domain density obtained by expanding $g(x)$ in terms
of the eigenfunctions $u(r,x)$ for the momentum operator considered above; doing so,
one obtains the density
\begin{align}
f_R(r)&=\left\vert G(r)\right\vert^2 =\frac{1}{2}\left(  \left\vert G_{1}%
(r)\right\vert ^{2}+\left\vert G_{2}(r)\right\vert ^{2} + 
 e^{i\beta}G_{1}^{\ast}(r)G_{2}(r)+e^{-i\beta}G_{1}(r)G_{2}^{\ast
}(r)\right) \nonumber
\end{align}
where
\[
G_{1,2}(r)=\frac{1}{\sqrt{2\pi }}\int_{-\infty}^\infty g_{1,2}(x)e^{-ixr}dx
\]
Note that, although $g_1(x)g_2(x)=0$, here we have $G_1(r)G_2(r) \ne 0$ because $G_1(r)$ and $G_2(r)$ are Fourier transforms of finite extent functions and therefore they extend over all $r$ \cite{slep2,marks}. Hence, unlike the case with $f_X(x)$,  the density $f_R(r)$ does depend on $\beta$. However, as Aharonov {\it et al.} showed,  
 the moments 
$\text{E\negthinspace}\left[ R^{n}\right]$  are
{\it independent} of $\beta$. 

Specifically, by virtue of the Fourier relations given above for
calculating moments, one has
\begin{align}
\text{E\negthinspace}\left[  R^{n}\right]  &  =\frac{1}{2}\int_{-\infty}^\infty\left(  g_{1}^{\ast
}(x)+e^{-i\beta}g_{2}^{\ast}(x)\right)  \left(  \frac{1}{i}\frac{d}%
{dx}\right)  ^{n}\left(  g_{1}(x)+e^{i\beta}g_{2}(x)\right)   dx  \nonumber \\
&  =\frac{1}{2}\int_{-\infty}^\infty g_{1}^{\ast}(x)\left(  \frac{1}{i}\frac{d}{dx}\right)
^{n}g_{1}(x)dx+\frac{1}{2}\int_{-\infty}^\infty g_{2}^{\ast}(x)\left(  \frac{1}{i}\frac{d}%
{dx}\right)  ^{n}g_{2}(x)dx  \nonumber\label{eq5}
\end{align}
where the final expression follows since $g_1(x)g_2(x)=0$.
 Hence, the density $f_R(r)$ is in general M-indeterminate. Historically, M-indeterminate densities were defined as those for which the moments do not uniquely determine the density, {\it and} for which all moments exist.  In their original considerations, Aharonov {\it et al.} did not address this latter criterion, but it has been addressed in the quantum context \cite{semon-taylor, semon-taylor-II, rafa2018, rafa2021} and we do so herein as well.

We generalize these ideas by considering other self-adjoint operators to obtain formulations for M-indeterminate densities. A key aspect of M-indeterminate densities that we take advantage of here is that, while M-indeterminate densities do not have a moment generating function,
 the characteristic function always exists.

\section{Characteristic function and self-adjoint operators} \label{prelim}

Let $R$ be a continuous real random variable with continuous probability density function
$f_R(r) \ (r\in \mathbb{R}^1) $.  
Then the characteristic function of  $f(r)$ is defined as \cite{feller,lukacs}
\begin{equation}
M({t})=\int_{-\infty}^{\infty}f_R(r)\,e^{i{t} r}\,dr \nonumber
\end{equation}
Given the characteristic function, one may obtain the probability density by Fourier inversion of the previous equation, namely
\begin{equation}
f_R(r)\,=\frac{1}{2\pi}\int_{-\infty}^{\infty}M({t})e^{-i{t} r}\,d{t} \nonumber
\end{equation}

\noindent\textbf{Theorem \ref{prelim}.1.} \emph{ If $f_R(r)=\left\vert G(r)\right\vert ^{2}$ and
$G(r)$ and $g(x)$ are related by Eqs. \eqref{eigen-expansion} and
\eqref{r-transform}, then the characteristic function may be obtained directly
in terms of $g(x)$ and the operator $\mathcal{A}$ as
\begin{equation}
M({t})=\int_{-\infty}^{\infty}g^{\ast}(x)\,e^{i{t}\mathcal{A}}\,g(x)\,dx
\label{general-char}%
\end{equation}
}
\noindent\textbf{Proof. }  See \cite{cohen1988, cohen2017} and references therein.  $\ $ \qedsymbol

\bigskip

\noindent\textbf{Corollary \ref{prelim}.1} \emph{For the particular case of the self-adjoint
operator
\begin{equation}
\mathcal{A}=\frac{1}{i}\frac{d}{dx} \nonumber
\end{equation}
we have that
\begin{equation}
M({t})=\int_{-\infty}^{\infty}g^{\ast}(x)e^{{t} \frac{d}{dx}%
}\,g(x)\,dx=\int_{-\infty}^{\infty}g^{\ast}(x)\,g(x+{t})\,dx \label{khinchin}
\end{equation}
}
\noindent{\bf Proof. }   The result follows readily by a Taylor series
expansion of $e^{{t}\frac{d}{dx}}$ \cite{wilcox}.  $\ $ \qedsymbol

\noindent{\bf Remark \ref{prelim}.1:} The right-hand side of Eq. \eqref{khinchin} is a well-known condition, first given by Khinchin \cite{khin,cohen1988}.

\section{A Stieltjes class characteristic function approach} \label{stieltjes-sec}
 One general formulation for generating so-called ``Stieltjes classess'' of M-indeterminate densities is
\cite{stoyanov2004,stoyanov2005}
\begin{equation}
f_{\varepsilon}(x)=f_{0}(x)\left[  1+{\varepsilon}\,h(x)\right]
\quad\label{Stieltjes}
\end{equation}
where $-1\leq{\varepsilon}\leq1$ is a real parameter, and $h(x)\neq0$
is a real, continuous, bounded ($|h(x)|\leq1$) function that satisfies the
constraint
\begin{equation}
\label{stieltjes-constraint}
\int_{-\infty}^{\infty}x^{n} f_{0}(x) h(x) dx = 0,
\ \ \ \ n\in I^{+}%
\end{equation}
by which it follows that the densities $f_\varepsilon(x)$  all have the same moments.

\subsection{Characteristic function constraint}
For the M-indeterminate densities of Eq. \eqref{Stieltjes},
we express the characteristic function as
\begin{align}
M_{\varepsilon}({t})  &  =\int_{-\infty}^{\infty}f_{\varepsilon
}(x)e^{i{t} x}dx\,=\, \int_{-\infty}^{\infty}f_{0}(x)e^{i{t}
x}dx+\varepsilon\int_{-\infty}^{\infty}h(x)f_{0}(x)e^{i{t} x}dx \nonumber \\
&  =M_{0}({t})\,+\,\varepsilon\,Q({t})  \nonumber \label{stieltjes-charfun}
\end{align}
where%
\begin{align}
M_{0}({t})  &  =\int_{-\infty}^{\infty}f_{0}(x)e^{i{t} x}dx \quad ; \quad
Q({t}) \,=\,  \int_{-\infty}^{\infty}f_{0}(x)h(x)e^{i{t} x}dx \nonumber 
\end{align}
Further, let $H({t})$ be the Fourier transform of $h(x)$,
\begin{equation}
{H({t})}=\int_{-\infty}^{\infty}h(x)e^{-i{t} x}dx \nonumber
\end{equation}

\bigskip

\noindent\textbf{Theorem \ref{stieltjes-sec}.1.} 
\emph{Let $M_{0}({t})$ and ${H({t})}$ and their derivatives 
$M_{0}^{(k)}= \frac{d^{k}}{d{t}^{k}}\,M_{0}(t)$ and $H^{(k)}= \frac{d^{k}}{d{t}^{k}}\,H(t)$ vanish
at $\pm\infty$.  Then for the probability densities $f_{\varepsilon}(x)$ of Eq. \eqref{Stieltjes},
the moments will be independent of $\varepsilon$ if
\begin{equation}
\label{charcondition-stieltjes}
\int_{-\infty}^{\infty}M_{0}^{(k)}
({t})\,H^{(n-k)}({t})\,d{t}=0,\quad k=0,1,...,n
\end{equation}
for all $n\in I^{+}$ and any value of $k$ as indicated.
}
\bigskip

\noindent {\bf Proof. }  The moments of $f_{\varepsilon}(x)$ are given by
\begin{equation}
\text{E\negthinspace}\left[  X^{n}\right]  \,=\,\int_{-\infty}^{\infty}%
x^{n}f_{\varepsilon}(x)dx\,=\,\int_{-\infty}^{\infty}x^{n}f_{0}%
(x)dx\,+\,\,\,\varepsilon\int_{-\infty}^{\infty}x^{n}f_{0}(x)h(x)dx \nonumber
\end{equation}
The moments can be equivalently obtained from the characteristic function by
\begin{equation}
\text{E\negthinspace}\left[  X^{n}\right]  \,=\,\left.  \left\{  
\frac{1}{i^n}\frac{d^{n}}{d{t}^{n}} \,M_{\varepsilon}({t})\right\}
\right\vert _{{t}=0}=\,\left.  \left\{   \frac{1}{i^{n}}\frac{d^{n}}%
{d{t}^{n}}\,M_{0}({t})\right\}  \right\vert _{{t}
=0}+\,\,\varepsilon\left.  \left\{   \frac{1}{i^{n}}\frac{d^{n}}{d{t}^{n}
}\,Q({t})\right\}  \right\vert _{{t}=0}%
\label{charfunmoments}%
\end{equation}
In order for the moments to be independent of $\varepsilon$, we require 
\begin{equation}
\left. \left\{  \frac{d^{n}}{d{t}^{n}}\,Q({t}) \right\} \right\vert _{{t}=0}=0 \nonumber
\end{equation}
Substituting in for $Q({t})$ and evaluating this expression, one can confirm that this condition is precisely that
given in Eq. \eqref{stieltjes-constraint}.  In terms of $M_0({t})$ and $H({t})$, straightforward evaluation yields
\begin{equation}
\left.  \left\{   \frac{d^{n}}{d{t}^{n}}\,Q({t}) \right\} \right\vert _{{t}
=0}\,=\,\left(  \frac{1}{2\pi}\right)  ^{2}\,\int_{-\infty}^{\infty}%
M_{0}^{(n)}({t}^{\prime}){H({t}^{\prime})}\,d{t}^{\prime}\,=\,0 \nonumber
\end{equation}
by which Eq. \eqref{charcondition-stieltjes} follows via integration by parts
as the characteristic function domain condition for the Stieltjes-classes of M-indeterminate densities  of Eq. \eqref{Stieltjes}.  $\ $ \qedsymbol

\subsection{Generating Stieltjes class densities via the characteristic
function}\label{charfun_approach}

The characteristic function constraint for M-indeterminate densities leads to
the following results:

\bigskip

\noindent\textbf{Theorem \ref{stieltjes-sec}.2.} 
{\it Let}
\begin{enumerate}
\item \emph{The characteristic function $M_{0}({t})$ be of finite extent:
$M_{0}({t}) = 0$ for $|{t}|>L>0$} \label{C1}
\item \emph{The function ${H({t})}=0$ for $|{t}|<L$, and is normalized such
that $\frac{1}{2\pi}\left|  \int_{-\infty}^{\infty}H({t})\,e^{i {t} x}\,
d{t}\right|  \,=\, |h(x)|\,\le\,1$} \label{C2}
\end{enumerate}
{\it Then the resulting densities $f_{\varepsilon}(x)$ of Eq. \eqref{Stieltjes} are independent of
$\varepsilon$.}

\bigskip
\noindent{\bf Corollary \ref{stieltjes-sec}.1:} \emph{ If all derivatives $\frac{d^n M_0({t})}{d{t}^n}$ are finite at ${t}=0$, then the densities $f_\varepsilon(x)$ are M-indeterminate.}

\bigskip
\noindent {\bf Proof. } 
The normalization on $H({t})$ insures that $f_{\varepsilon}(x)\geq0$. Then,
imposing the finite support constraint on $M_{0}({t})$ yields, via Eq.
\eqref{charcondition-stieltjes} with $k=0$,
\begin{equation}
\int_{-\infty}^{\infty}M_{0}({t}^{\prime})\,{H^{(n)}({t}^{\prime}%
)}\,d{t}^{\prime}\,=\,\int_{-L}^{L}M_{0}({t}^{\prime})\,{H^{(n)}%
({t}^{\prime})}\,d{t}^{\prime}  \nonumber
\end{equation}
Since $H({t}^{\prime})=0$ over the limits of integration, the integral
equals zero;
hence, the moments are independent of $\varepsilon$ but the densities
$f_{\varepsilon}(x)$ are not.  Finally, the moments may be obtained by
$ \left. \left\{  \frac{1}{i^n}\frac{d^{n}}{d{t}^{n}}\,M_\varepsilon({t}) \right\} \right\vert _{{t}=0}$=
$ \left. \left\{  \frac{1}{i^n}\frac{d^{n}}{d{t}^{n}}\,M_0({t}) \right\} \right\vert _{{t}=0}. $ 
By virtue of the corollary, all moments are finite.  Hence, the densities are M-indeterminate.
 $\ $ \qedsymbol

\bigskip

\noindent These results point to a simple procedure for generating an unlimited number of
Stieltjes classes of  M-indeterminate densities:

\bigskip

Choose \textit{any} infinitely differentiable, finite-extent function $g(x)$, normalized
to 1; namely,
\begin{align}
\frac{d^{n}g(x)}{dx^{n}}  &  <\infty,\,\,\,\,n=0,1,2,... \nonumber \\
g(x)  &  =0,\qquad|x|>\frac{L}{2}>0 \nonumber \\
\int_{-\infty}^{\infty}|g(x)|^{2}\,dx  &  =1  \nonumber
\end{align}
Examples of such functions are ``bump functions'' \cite{tu}.
Then calculate
\begin{equation}
\label{ordinary-charfun}
M_{0}({t})=\int_{-\infty}^{\infty}g^{\ast
}(x)g(x+{t})\ dx \nonumber
\end{equation}

The density function is
\begin{align}
f_{0}(x) &  =\frac{1}{2\pi}\int_{-\infty}^{\infty}M_{0}({t})\,e^{-i{t}
x}\,d{t}\label{density_0} \nonumber \\
&  =\frac{1}{2\pi}\int_{-\infty}^{\infty}\int_{-\infty}^{\infty}g^{\ast
}(y)g(y+{t})\,dy\,e^{-i{t} x}\,d{t} \nonumber \\
&=\frac{1}{{2\pi}}\left|\int_{-\infty}^{\infty}g(y)e^{-ixy}\,dy \right|^2 \nonumber
\end{align}
Note that because $g(x)$ is normalized, $f_{0}(x)$ is likewise
normalized, {\it i.e.,} 
$\int_{-\infty}^{\infty}f_{0}(x)\,dx\,=\,
\int_{-\infty}^{\infty}|g(x)|^{2}\,dx\,=\,1.   
$
Note further that, by virtue of the finite
extent of $g(x)$, the characteristic function $M_0(t)$ is also of finite extent and in
particular it equals zero for $|{t}|>L$.  Hence the density $f_{0}(x)$ always exists.

Then, form the Stieltjes classes of M-indeterminate densities
\begin{equation}
f_{\varepsilon}(x;\,\lambda,\phi)=f_{0}(x)\left[  1+\varepsilon\,\cos\left(
\lambda x+\phi\right)  \right]  ,\quad\lambda>L,\quad-1\leq\varepsilon
\leq1,\quad-\pi\leq\phi\leq\pi\label{Method1-Stieltjes-class} 
\end{equation}
Because $\lambda>L$ and the characteristic function equals zero for
$|{t}|>L$, it follows that $\int_{-\infty}^{\infty}x^{n}f_{0}%
(x)\cos\left(  \lambda x+\phi\right)  dx=0$ (which is more readily evaluated
in the characteristic function domain via Eq. \eqref{charfunmoments}). Hence,
$f_{0}(x)$ and the densities $f_{\varepsilon}(x;\,\lambda,\phi)$ all
have identical moments E$[X^{n}]$, consistent with Theorem \ref{stieltjes-sec}.2.
\bigskip

\noindent{\bf Remark \ref{stieltjes-sec}.1:} Because $g(x)$ is infinitely differentiable, 
it follows that all derivatives $\frac{d^n M({t})}{d{t}^n}$ are finite at ${t}=0$, and hence all moments exist, which is a necessary condition for a density to be classified as M-indeterminate.

\noindent{\bf Remark \ref{stieltjes-sec}.2:}  This form of the characteristic function corresponds to the special case considered in Corollary \ref{prelim}.1.

\section{Generalized characteristic function and operator approach} \label{generalized}
We now extend the characteristic function approach by making use of Eq.
\eqref{general-char} with general self-adjoint operators ${\cal{A}}$ and taking the
(normalized) function to be
\begin{equation}
g(x;\,\beta)=g_{1}(x)+e^{i\beta}g_{2}(x), \quad \hbox{with} \quad g_{1}(x)g_{2}(x)=0,
 \quad  \int_{-\infty}^\infty |g_{1}(x)|^2dx=\int_{-\infty}^\infty  |g_{2}(x)|^2dx=\frac{1}{2} \label{m-ind-g}%
\end{equation}
where we have made the dependence on the real parameter  $\beta$ explicit in our notation.
For clarity, we shall denote here the characteristic function of Eq.
\eqref{general-char} by $M_{\mathcal{A}}({t};\, \beta)$, with corresponding probability
density $f_{\mathcal{A}}(r;\,\beta)$ for the continuous real random variable $R$.

\noindent \textbf{Theorem \ref{generalized}.1. }
\emph{Let $g(x;\,\beta)$ be given by Eq. \eqref{m-ind-g} 
and let the characteristic function
$M_{\mathcal{A}}({t};\,\beta)$ be given by Eq. \eqref{general-char}.   
Further, let the operator ${\mathcal{A}}$ be such that the support of
$\mathcal{A}^{n}g(x;\,\beta)$ is the same as that of $g(x;\,\beta)$.
Then, the moments $\text{E\negthinspace}\left[  R^{n}\right] $ of the family of
densities%
\begin{align}
f_{\mathcal{A}}(r;\,\beta)  &  =\frac{1}{2\pi}\int_{-\infty}^{\infty}M_{\mathcal{A}%
}({t};\,\beta)\,e^{-i{t} r}\,d{t} \nonumber \\
&  =\frac{1}{2\pi}\int_{-\infty}^{\infty}\int_{-\infty}^{\infty}g^{\ast
}(x;\,\beta)e^{i{t}\mathcal{A}}\,g(x;\,\beta)\,dx\,e^{-i{t} r}\,d{t} \label{gen-m-ind-densities}
\end{align}
are independent of the parameter $\beta$.
}

\noindent {\bf Proof. } 
The
moments are
\begin{align}
\text{E\negthinspace}\left[  R^{n}\right]  &=  \int_{-\infty}^\infty r^n\,f_{\mathcal{A}}(r;\,\beta)\,dr
\,=\, \left.  \left\{  \left(
\frac{1}{i}\frac{d}{d{t}}\right)  ^{n}\,M_{\mathcal{A}}({t};\,\beta)\right\}
\right\vert _{{t}=0}\,  \nonumber \\
&= \int_{-\infty}^{\infty}g_{1}^{\ast
}(x)\,\mathcal{A}^{n}\,g_{1}(x)\,dx\,+\,\int_{-\infty}^{\infty}g_{2}^{\ast
}(x)\,\mathcal{A}^{n}\,g_{2}(x)\,dx \nonumber \\
& +\,e^{i\beta}\, \int_{-\infty}^{\infty}g_{1}^{\ast}(x)\,\mathcal{A}^{n}%
\,g_{2}(x)\,dx\,+\,e^{-i\beta}\, \int_{-\infty}^{\infty}g_{2}^{\ast}(x)\,
\mathcal{A}^{n}\,g_{1}(x)\,dx \label{method3_moments}  \nonumber
\end{align}
But since $g_{1}(x)g_{2}(x)=0$ and $\mathcal{A}^{n}g_{1}(x)$ has the same support
as $g_{1}(x)$, and likewise for $g_2(x)$, it follows that 
\begin{equation}
 e^{i\beta}\, \int_{-\infty}^{\infty}g_{1}^{\ast}(x)\,\mathcal{A}^{n}%
\,g_{2}(x)\,dx\,+\,e^{-i\beta}\, \int_{-\infty}^{\infty}g_{2}^{\ast}(x)\,
\mathcal{A}^{n}\,g_{1}(x)\,dx \,=\,0 \nonumber
\end{equation}
Accordingly,
the moments are independent of $\beta$.    $\ $ \qedsymbol

\noindent \textbf{Corollary \ref{generalized}.1. }
\emph{In addition to the conditions of Theorem \ref{generalized}.1, let
\begin{align}
m_1^{(n)} &=  \int_{-\infty}^{\infty}g_{1}^{\ast}(x)\,\mathcal{A}^{n}\,g_{1}(x)\,dx < \infty \nonumber \\
m_2^{(n)} &= \int_{-\infty}^{\infty}g_{2}^{\ast}(x)\,\mathcal{A}^{n}\,g_{2}(x)\,dx < \infty \nonumber
\end{align}
for all $n \in I^+$. 
Further, let the support of
$e^{i{t}\mathcal{A}}\,g(x)$ be different than the support of $g(x)$ (e.g., recall Corollary \ref{prelim}.1).
Then the densities defined by Eq. \eqref{gen-m-ind-densities} are M-indeterminate.
}

\noindent {\bf Proof.}
Theorem \ref{generalized}.1 proves that the densities $f_{\mathcal{A}}(r;\,\beta)$ have identical moments,
which is a necessary condition to be M-indeterminate.  A second necessary condition is that all moments exist; the conditions on $m_1^{(n)}$ and $m_2^{(n)}$ ensure that all moments are finite.  What remains to show is that $f_{\mathcal{A}}(r;\,\beta)$ (equivalently, $M_{\mathcal{A}}({t};\,\beta) $) depends on $\beta$.  
Substituting \eqref{m-ind-g} into \eqref{general-char} yields
\begin{align}
M_{\mathcal{A}}({t};\,\beta)  &  =\int_{-\infty}^{\infty}g_{1}^{\ast}%
(x)\,e^{i{t}\mathcal{A}}\,g_{1}(x)\,dx\,+\,\int_{-\infty}^{\infty}%
g_{2}^{\ast}(x)\,e^{i{t}\mathcal{A}}\,g_{2}(x)\,dx\, \nonumber \\
&  +\,e^{i\beta}\,\int_{-\infty}^{\infty}g_{1}^{\ast}(x)\,e^{i{t}
\mathcal{A}}\,g_{2}(x)\,dx\,+\,e^{-i\beta}\,\int_{-\infty}^{\infty}g_{2}^{\ast
}(x)\,e^{i{t}\mathcal{A}}\,g_{1}(x)\,dx  \label{general-charfun} \nonumber
\end{align}
Now, even though $g_1(x)g_2(x)=0$, by the conditions of the theorem 
it follows that $g_{1}^{\ast}(x)\,e^{i{t}\mathcal{A}}\,g_{2}(x) \ne 0$ and likewise $g_{2}^{\ast}(x)\,e^{i{t}\mathcal{A}}\,g_{1}(x) \ne 0$,
such that the last two integrals are nonzero.  Hence,  in general the densities  depend on $\beta$ but
the moments do not and are finite, and therefore the densities $f_{\mathcal{A}}(r;\,\beta)$ are M-indeterminate.  $\ $ \qedsymbol

\noindent \textbf{Remark \ref{generalized}.1: } All differential operators $\mathcal{A}$ and infinitely differentiable functions $g(x;\,\beta)$ given by Eq. \eqref{m-ind-g} 
satisfy the conditions of Corollary  \ref{generalized}.1
\cite{rafa2021}.

\subsection{Example 1}
For the case where $\mathcal{A}=\frac{1}{i}\frac
{d}{dx}$, the characteristic function is (see Corollary \ref{prelim}.1)
\begin{align}
M_{\mathcal{A}}({t};\,\beta) &  =\int_{-\infty}^{\infty}g_{1}^{\ast}(x)\,g_{1}%
(x+{t})\,dx\,+\,\int_{-\infty}^{\infty}g_{2}^{\ast}(x)\,g_{2}%
(x+{t})\,dx \nonumber \\
&  +\,e^{i\beta}\,\int_{-\infty}^{\infty}g_{1}^{\ast}(x)\,g_{2}(x+{t}
)\,dx\,+\,e^{-i\beta}\,\int_{-\infty}^{\infty}g_{2}^{\ast}(x)\,g_{1}%
(x+{t})\,dx  \nonumber
\end{align}
which clearly depends on $\beta$ since the latter two integrals are nonzero, even though 
$g_{1}(x)g_{2}(x)=0.$
The moments are
\begin{align}
\text{E\negthinspace}\left[  R^{n}\right]  \,=\,\left.  \left\{  \left(
\frac{1}{i}\frac{d}{d{t}}\right)  ^{n}\,M_{\mathcal{A}}({t};\,\beta)\right\}
\right\vert _{{t}=0}\, &  =\,\,\,\int_{-\infty}^{\infty}g_{1}^{\ast
}(x)\,\left(  \frac{1}{i}\frac{d}{dx}\right)  ^{n}g_{1}(x)\,dy\,+\,\int
_{-\infty}^{\infty}g_{2}^{\ast}(x)\,\left(  \frac{1}{i}\frac{d}{dx}\right)
^{n}g_{2}(x)\,dy\,+\, \nonumber \\
&  e^{i\beta}\,\int_{-\infty}^{\infty}g_{1}^{\ast}(x)\,\left(  \frac{1}%
{i}\frac{d}{dx}\right)  ^{n}g_{2}(x)\,dy\,+\,e^{-i\beta}\,\int_{-\infty}^{\infty}%
g_{2}^{\ast}(x)\,\left(  \frac{1}{i}\frac{d}{dx}\right)
^{n}g_{1}(x)\,dx \nonumber
\end{align}
which are independent of $\beta$ since $g_{1}(x)g_{2}(x)=0$.

\noindent \textbf{Remark \ref{generalized}.2: }The operator $\mathcal{A}=\frac{1}{i}\frac
{d}{dx}$ is fundamental in quantum mechanics as it is associated with momentum.

\subsection{Example 2}
The operator 
\begin{equation}
\mathcal{A}=\frac{1}{i}\frac{d}{dx}+c(n+1)x^{n} \nonumber
\end{equation}
where $c\in R$ and $n\in I^+$, satisfies the conditions of Theorem 4.1 and Corollary 4.1; namely,
it is self-adjoint, and the support of $\mathcal{A}^ng(x)$ equals that of $g(x)$,
whereas the support of $e^{i {t} \mathcal{A}}g(x)$ differs from the support of $g(x)$.

For $g(x)$ as in Eq. \eqref{m-ind-g}, we have by Eq. \eqref{general-char}
that the characteristic function is
\begin{align}
M_{\mathcal{A}}({t};\,\beta) 
&  =\int_{-\infty}^\infty\left(  g_{1}^{\ast}(x)+e^{-i\beta}g_{2}
^{\ast}(x)\right)  e^{i{t}\mathcal{A}}\left(  g_{1}(x)+e^{i\beta}f
_{2}(x)\right)  dx \nonumber\\
&  = M_{11}({t})+M_{22}({t})+e^{i\beta}M_{12}%
({t})+e^{-i\beta}M_{21}({t}) \nonumber
\end{align}
where
\begin{align}
M_{lm}({t}) &  =\int_{-\infty}^\infty g_{l}^{\ast}(x)e^{i{t}\mathcal{A}}g_{m}(x)dx \nonumber
\end{align}
Note that the characteristic function, and hence the density, depends on $\beta$.
For the moments to be independent of $\beta$, 
we require that
\begin{equation} \label{special2-condition}
 \left. \left\{e^{i\beta}M_{12}^{(n)}({t})+e^{-i\beta}M_{21}%
^{(n)}({t})\right\} \right\vert_{{t}=0}=0  \nonumber
\end{equation}
where
\begin{equation}
M_{lm}^{(n)}({t})\,=\,\frac{d^n}{d {t}^n}M_{lm}({t}) \nonumber
\end{equation}
Evaluating this expression, we have
\begin{align}
\left\{ e^{i\beta}M_{12}^{(n)}({t})+e^{-i\beta}M_{21}^{(n)}({t})\right\}
_{{t}=0} &  =\left\{  e^{i\beta}\int_{-\infty}^\infty g_{1}^{\ast}(x)i^{n}\mathcal{A}%
^{n}e^{i{t}\mathcal{A}}g_{2}(x)dx+e^{-i\beta}\int_{-\infty}^\infty g_{2}^{\ast}%
(x)i^{n}\mathcal{A}^{n}e^{i{t}\mathcal{A}}g_{1}(x)dx\right\}  _{{t}
=0} \nonumber \\
&  =e^{i\beta}i^{n}\int_{-\infty}^\infty g_{1}^{\ast}(x)\mathcal{A}^{n}g_{2}%
(x)dx+e^{-i\beta}i^{n}\int_{-\infty}^\infty g_{2}^{\ast}(x)\mathcal{A}^{n}g_{1}(x)dx \nonumber \\
&=0 \nonumber
\end{align}
where the last step follows since $g_1(x)g_2(x)=0$ and  the support of $\mathcal{A}^ng_1(x)$ equals that of $g_1(x)$, and similarly for $\mathcal{A}^ng_2(x)$. 

\noindent \textbf{Remark \ref{generalized}.3: }Special cases of this operator  are the creation and annihilation operators for the quantum harmonic oscillator.

\section{Closing remarks}
We have developed new methods to construct M-indeterminate probability densities via the characteristic function together with operator methods used in quantum mechanics.
We now briefly remark on applications as appropriate for this journal. Generally speaking, the M-indeterminate moment problem has been formulated with the necessary condition that all moments exist.  Of course, there are important densities for which all moments do not exist. The procedures we have given for generating an unlimited number of densities that depend on a parameter, whereas the moments do not, is applicable in general, regardless of whether all moments exist (as of course is the approach of Eq. \eqref{Stieltjes}). 

Since the density and the characteristic function are Fourier transform pairs, the M-indeterminate condition that all moments exist implies that the characteristic function is infinitely differentiable at the origin.   Using so-called “bump” functions, as noted in section \ref{charfun_approach}, ensures that this condition is met.  The existence of all moments coupled with the Fourier relation between the density and the characteristic function also has implications on the tails of the density.  Practical issues involving the number of moments actually known and possible errors in the moments are important questions.  Also, the general issue of the differentiability of a function from a practical point of view, and relatedly the estimation of the tails of a density as well as the influence of the tails on physical properties, has been addressed by a number of authors (e.g. \cite{heyde,slep1,slep2}).


\begin{thebibliography}{99}                                                                                               %

\small{

\bibitem {ahar1}Aharonov, Y., Pendleton, H. and Petersen, A., 1969. Modular
Variables in Quantum Theory. Int. J. Th. Phys. 2(3), 213-230.

\bibitem {ahar2}Aharonov, Y., Pendleton, H. and Petersen, A., 1970.
Deterministic Quantum Interference
Experiments. Int. J. Th. Phys. 3(6), 443-442.

\bibitem {ahar-book}Aharonov, Y. and Rohrlich, D., 2005. Quantum Paradoxes.
Wiley-VCH.

\bibitem {bohm}Bohm, D., 1951. Quantum Theory. Prentice-Hall, New York.

\bibitem {cohen1988}Cohen, L., 1988. Rules of probability in quantum
mechanics. Found. Phys. 18(10), 983-998.

\bibitem {cohen2017}Cohen, L., 2017. Are there quantum operators and wave
functions in standard probability theory?, in: Wong, M. W., Zhu, H. (Eds.),
Pseudo-Differential Operators: Groups, Geometry and Applications, pp. 133-147.
Birkh\"{a}user Mathematics.

\bibitem {feller}Feller, W., 1971. An Introduction to Probability Theory and
Its Applications, Vol. 2. John Wiley and Sons, New York.

\bibitem{heyde} Heyde, C. and Kou, S., 2004. On the controversy over tailweight of distributions. Oper. Res. Letters 32, 399-408.

\bibitem{khin}Khinchin, A., 1937. On a property of characteristic functions. 
Bull. Univ. Moscow, Vol. 1.

\bibitem {lukacs}Lukacs, E., 1970. Characteristic Functions, 2nd ed. Griffin
\& Co., London.

\bibitem {marks}Marks II, R. J., 2009. Handbook of Fourier Analysis \& Its
Applications. Oxford Univ. Press.

\bibitem {merz}Merzbacher, E., 1998.  Quantum Mechancs. John Wiley \& Sons, Inc.

\bibitem {morse}Morse, P. and Feshbach, H., 1953. Methods of Theoretical
Physics, Part I. McGraw-Hill Book Company, New York.

\bibitem {rafa2018} Sala Mayato, R., Loughlin, P. and Cohen, L., 2018.
M-indeterminate distributions in quantum mechanics and the non-overlapping
wave function paradox. Physics Letters A 382, 2914--2921.

\bibitem {rafa2021} Sala Mayato, R., Loughlin, P. and Cohen, L., 2022. Generating
M-indeterminate probability densities by way of quantum mechanics. J. Theor. Probab., 35, 1537-1555.

\bibitem {semon-taylor} 
Semon, M.D., Taylor, J.R., 1987. Expectation values in the Aharonov-Bohm effect.
Il Nuovo Cimento B 97(1), 25-40.

\bibitem {semon-taylor-II} 
Semon, M.D. and Taylor, J.R., 1987. Expectation values in the Aharonov-Bohm effect. - II.
Il Nuovo Cimento B 100(3), 389-401.

\bibitem{slep1} Slepian, D., 1976. On bandwidth. Proc. IEEE, 64(3), 292-300.

\bibitem{slep2} Slepian, D., 1983. Some comments on Fourier analysis, uncertainty and modeling. SIAM Review 25(3), 379-393.

\bibitem {stoyanov2004}Stoyanov, J., 2004. Stieltjes classes for
moment-indeterminate probability distributions. J. Applied Probability 41A, 281-294.

\bibitem {stoyanov2005}Stoyanov, J. and Tolmatz, L., 2005. Method for
constructing Stieltjes classes for m-indeterminate probability distributions.
Appl. Math. and Comp. 165, 669-685.

\bibitem {stoyanov2013}Stoyanov, J., 2013. Counterexamples in Probability, 3rd
ed. Dover.


\bibitem {tu}Tu, L. W., 2011. An Introduction to Manifolds, Second Edition.
Springer, NY.


\bibitem {wilcox}Wilcox, R.M., 1967. Exponential operators and parameter
differentiation in quantum physics. J. Math. Phys. 8, 962-981.


} 
\end{thebibliography}
\end{document}